\documentclass[reqno]{amsart}
\usepackage{hyperref}

\begin{document}

\title[\hfilneg \hfil Uniqueness on Meromorphic function concerning their differential-difference operators]
{Uniqueness on Meromorphic function concerning their differential-difference operators}

\author[XiaoHuang Huang \hfil \hfilneg]
{XiaoHuang Huang}

\address{XiaoHuang Huang: Corresponding author\newline
Department of Mathematics, Southern University of Science and Technology, Shenzhen 518055, China}
\email{1838394005@qq.com}

\begin{abstract}
In this paper, we study the uniqueness of the differential-difference of meromorphic functions. We prove the following result: Let $f$ be a  nonconstant meromorphic function of $\rho_{2}(f)<1$, let $\eta$ be a non-zero complex number, $n\geq1, k\geq0$  two  integers and let $a\not\equiv0,\infty$  be a small function of $f$. If $f$ and $(\Delta_{\eta}^{n}f)^{(k)}$ share $0,\infty$ CM and share $a$ IM, then $f\equiv(\Delta_{\eta}^{n}f)^{(k)}$, which use a completely different method to  improve some results due to Chen-Xu \cite{cx1}.
\end{abstract}

\maketitle
\numberwithin{equation}{section}
\newtheorem{theorem}{Theorem}[section]
\newtheorem{lemma}[theorem]{Lemma}
\newtheorem{remark}[theorem]{Remark}
\newtheorem{corollary}[theorem]{Corollary}
\newtheorem{example}[theorem]{Example}
\newtheorem{problem}[theorem]{Problem}
\allowdisplaybreaks

\section{Introduction and main results}
In this paper, we use the standard denotations in the  Nevanlinna value distribution theory, see([9,17,18]).
 Throughout this paper,  $f(z)$ is a meromorphic function on the whole complex plane. $S(r,f)$ means that $S(r, f) = o(T(r, f))$, as $r\to \infty $ outside of a possible exceptional set of finite logarithmic measure.
 Define
 $$\rho(f)=\varlimsup_{r\rightarrow\infty}\frac{log^{+}T(r,f)}{logr},$$
 $$\rho_{2}(f)=\varlimsup_{r\rightarrow\infty}\frac{log^{+}log^{+}T(r,f)}{logr}$$
by the order, lower order  and the hyper-order  of $f$, respectively.

A meromorphic function $a(z)$ satisfying $T(r,a)=S(r,f)$ is called a small function of $f$. We say that two nonconstant meromorphic functions $f$ and $g$ share small function $a$ CM(IM) if $f-a$ and $g-a$ have the same zeros counting multiplicities (ignoring multiplicities).

Let $f(z)$ be a meromorphic function, and a finite complex number $\eta$. we define  its difference operators by
\begin{equation*}
\Delta_\eta f(z)=f(z+\eta)-f(z), \quad \Delta_\eta^{n}f(z)=\Delta_{\eta}^{n-1}(\Delta_\eta f(z)).
\end{equation*}

In recent years, there has been tremendous interests in developing  the value distribution of meromorphic functions with respect to difference analogue.

In 2011, Heittokangas et al [10] proved a uniqueness Theorem of meromorphic function concerning shift.

\

{\bf Theorem A}
 Let $f(z)$ be a nonconstant meromorphic function of finite order, let $\eta$ be a nonzero finite complex value, and let $a, b, c$ be three  distinct complex values in the extended plane. If $f(z)$ and $f(z+\eta)$ share $a, b, c$ CM, then $f(z)\equiv f(z+\eta).$

 Corresponding to Theorem 3, Cui-Chen [5], Li-Yi-Kang [11], L\"u-L\"u [13] proved

\

{\bf Theorem B}
 Let $f(z)$ be a transcendental meromorphic function of finite order, let $\eta$ be a non-zero complex number, and let $a$ and $b$ be two distinct complex numbers. If $f(z)$ and $\Delta_{\eta}f(z)$ share $a$, $b$, $\infty$ CM, then $f(z)\equiv \Delta_{\eta}f(z)$.

It's naturally to ask a question from  Theorem B that

\

{\bf Question 1}
 Can "3 CM" be replaced by "2 CM+1 IM"?

In 2020, Chen-Xu [1] proved.

\

{\bf Theorem C}
 Let $f$ be a  transcendental entire function of $\rho_{2}(f)<1$, let $\eta$ be a non-zero complex number, $n\geq1$ a integer and let $a\not\equiv0,\infty$  be a complex value. If $f$ and $\Delta_{\eta}f$ share $0,\infty$ CM and share $a$ IM, then $f\equiv\Delta_{\eta}f$.

In real analysis, the time-delay differential equation $f'(x)=f(x-k), k>0$ has been extensively studied. As for a complex variable counterpart, Liu and Dong
[12]studied the complex differential-difference equation $f'(z)=f(z+c)$,  where c is a non-zero constant. In [15], Qi-Li-Yang investigated the value
sharing problem related to $f'(z)$ and $f(z+c)$, and proved

\

{\bf Theorem D}
 Let $f$ be a nonconstant entire function of finite order, and let $a, c$ be two nonzero finite complex values.
If $f'(z)$ and $f(z+c)$ share $0, a$ CM, then $f'(z)\equiv f(z+c).$

Recently, Qi and Yang [16] improved Theorem E and proved

\

{\bf Theorem E}
 Let $f$ be a nonconstant entire function of finite order, and let $a, c$ be two nonzero finite complex values.
If $f'(z)$ and $f(z+c)$ share $0$ CM and $a$ IM, then $f'(z)\equiv f(z+c).$

According to Theorem C and Theorem E, we can pose a question.
\

{\bf Question 2}
 Let $f$ be a  nonconstant meromorphic function of $\rho_{2}(f)<1$, let $\eta$ be a non-zero complex number, $n\geq1, k\geq0$  two  integers and let $a\not\equiv0,\infty$  be a small function of $f$. If $f$ and $(\Delta_{\eta}^{n}f)^{(k)}$ share $0,\infty$ CM and share $a$ IM, is $f(z)\equiv(\Delta_{\eta}^{n}f(z))^{(k)}$?

We give a positive answer to above question. We prove.

\

{\bf Theorem 1}
 Let $f$ be a   transcendental meromorphic function of $\rho_{2}(f)<1$, let $\eta$ be a non-zero complex number, $n\geq1, k\geq0$  two  integers and let $a\not\equiv0,\infty$  be a rational function. If $f$ and $(\Delta_{\eta}^{n}f)^{(k)}$ share $0,\infty$ CM and share $a$ IM, then $f\equiv(\Delta_{\eta}^{n}f)^{(k)}$.

we will see in section 3, there is nothing to do with $k$, that is to say when $k=0$, we can get the following corollary.

\

{\bf Corollary }
 Let $f$ be a   transcendental meromorphic function of $\rho_{2}(f)<1$, let $\eta$ be a non-zero complex number, $n\geq1$ a integer and let $a\not\equiv0,\infty$  be a rational function. If $f$ and $\Delta_{\eta}^{n}f$ share $0,\infty$ CM and share $a$ IM, then $f\equiv\Delta_{\eta}^{n}f$.

\

{\bf Example 1} Let $f(z)=e^{iz}$, and $\eta=\frac{\pi}{2}$. Then $(\Delta_{\eta}^{n}f(z))^{(k)}=i^{n}(i-1)^{n}e^{iz}$. Easy to see that  $f(z)$ and $(\Delta_{\eta}^{n}f(z))^{(k)}$ share $0$ and $\infty$ CM, but $f(z)\not\equiv(\Delta_{\eta}^{n}f(z))^{(k)}$.

\

{\bf Example 2} Let $f(z)=e^{z}$, and $\eta=log4$. Then $(\Delta_{\eta}^{n}f(z))^{(k)}=3^{n}e^{z}$. Easy to see that  $f(z)$ and $(\Delta_{\eta}^{n}f(z))^{(k)}$ share $0$ and $\infty$ CM, but $f(z)\not\equiv(\Delta_{\eta}^{n}f(z))^{(k)}$.

\

{\bf Example 3} Let $f(z)=\frac{e^{z}}{e^{2z}-1}$, and $\eta=i\pi$. Then $\Delta_{\eta}^{n}f(z)=\frac{(-2)^{n}e^{z}}{e^{2z}-1}$. Easy to see that  $f(z)$ and $\Delta_{\eta}^{n}f(z)$ share $0$ and $\infty$ CM, but $f(z)\not\equiv\Delta_{\eta}^{n}f(z)$.

Above three  examples show that the  conditions "2 CM+1 IM" is necessary.

A natural question from Theorem 1 is that what happens if sharing $0$ CM is replaced by sharing a nonzero value $a$ CM? Hence we raise a question.

\

{\bf Question} Let $f$ be a  transcendental meromorphic function of $\rho_{2}(f)<1$, let $\eta$ be a non-zero complex number, $n\geq1, k\geq0$  two  integers and let $a\neq0$ be a finite complex number, and  $b\not\equiv\infty$  a rational function. If $f$ and $(\Delta_{\eta}^{n}f)^{(k)}$ share $a,\infty$ CM and share $b$ IM, is  $f\equiv(\Delta_{\eta}^{n}f)^{(k)}$ still valid?

\section{Some Lemmas}

\begin{lemma}\label{21l}\cite{h3} Let $f$ be a nonconstant meromorphic function of $\rho_{2}(f)<1$,  and let $\eta$ be a non-zero complex number. Then
$$m(r,\frac{f(z+\eta)}{f(z)})=S(r, f),$$
for all r outside of a possible exceptional set E with finite logarithmic measure.
\end{lemma}

\begin{lemma}\label{21l}\cite{h3}
Let $f$ be a nonconstant meromorphic function of $\rho_{2}(f)<1$, and let $\eta\neq0$ be a finite complex number. Then
$$T(r,f(z+\eta))=T(r, f(z))+S(r,f).$$
\end{lemma}

\begin{lemma}\label{23l}\cite{h4,y1}
Let $f$ and $g$ be two nonconstant meromorphic functions with $\rho(f)$ and $\rho(g)$ as their orders respectively. Then $$\rho(fg)\leq max\{\rho(f),\rho(g)\}.$$
\end{lemma}

\begin{lemma}\label{24l}\cite{h4,y1}
Let $f=e^{h}$  be an entire functions with $\rho_{2}(f)<1$, where $h$ is a polynomial. Then $\rho(f)=deg(h)$, and $\rho(f)=\lambda(f)$.
\end{lemma}

\section{The proof of Theorem 1 }
 Assume that $f\not\equiv(\Delta_{\eta}^{n}f)^{(k)}$. Since $f$ is a transcendental meromorphic function of $\rho_{2}(f)<1$, $f$ and $(\Delta_{\eta}^{n}f)^{(k)}$ share $0,\infty$ CM, then there is a nonzero polynomial $p$  such that
\begin{eqnarray}
\frac{(\Delta_{\eta}^{n}f)^{(k)}}{f}=e^{p},
\end{eqnarray}
then by Lemma 2.1
\begin{eqnarray}
T(r,e^{p})=m(r,e^{p})=m(r,\frac{(\Delta_{\eta}^{n}f)^{(k)}}{f})=S(r,f).
\end{eqnarray}
On the other hand, (3.1) can be rewritten as
\begin{eqnarray}
\frac{(\Delta_{\eta}^{n}f)^{(k)}-f}{f}=e^{p}-1,
\end{eqnarray}
then from the fact that $f$ and $(\Delta_{\eta}^{n}f)^{(k)}$ share $a$ IM, we get
\begin{eqnarray}
\overline{N}(r,\frac{1}{f-a})\leq N(r,\frac{1}{e^{p}-1})\leq T(r,e^{p})=S(r,f).
\end{eqnarray}
From the fact that $f$ and $(\Delta_{\eta}^{n}f)^{(k)}$ share $0$ CM, then one has
\begin{eqnarray*}
\begin{aligned}
m(r,\frac{1}{f})+m(r,\frac{1}{f-a})&\leq m(r,\frac{1}{(\Delta_{\eta}^{n}f)^{(k)}})+m(r,\frac{1}{(\Delta_{\eta}^{n}(f-a))^{(k)}})+S(r,f)\\
&\leq 2T(r,(\Delta_{\eta}^{n}f)^{(k)})-N(r,\frac{1}{(\Delta_{\eta}^{n}(f-a))^{(k)}})\\
&-N(r,\frac{1}{f})+S(r,f),
\end{aligned}
\end{eqnarray*}
which implies
\begin{eqnarray}
2T(r,f)\leq 2T(r,(\Delta_{\eta}^{n}f)^{(k)})-N(r,\frac{1}{(\Delta_{\eta}^{n}(f-a))^{(k)}})+S(r,f).
\end{eqnarray}
We can deduce from (3.1) that
\begin{eqnarray}
T(r,f)= T(r,(\Delta_{\eta}^{n}f)^{(k)})+S(r,f).
\end{eqnarray}
Combing (3.5) with (3.6), we get
\begin{eqnarray}
N(r,\frac{1}{(\Delta_{\eta}^{n}(f-a))^{(k)}})= S(r,f).
\end{eqnarray}
Set
\begin{eqnarray}
\psi=\frac{(\Delta_{\eta}^{n}(f-a))^{(k)}}{f-a}.
\end{eqnarray}
Easy to see that
\begin{eqnarray}
N(r,\frac{1}{\psi})\leq N(r,\frac{1}{(\Delta_{\eta}^{n}(f-a))^{(k)}})+N(r,\frac{1}{a})= S(r,f),
\end{eqnarray}
\begin{eqnarray}
N(r,\psi)\leq N(r,\frac{1}{f-a})+N(r,(\Delta_{\eta}^{n}a)^{(k)})= S(r,f).
\end{eqnarray}
Hence by Lemma 2.2, we have
\begin{align}
T(r,\psi)&=m(r,\psi)+N(r,\psi)\notag\\
&\leq m(r,\frac{(\Delta_{\eta}^{n}(f-a))^{(k)}}{f-a})+N(r,\frac{1}{f-a})+S(r,f)\notag\\
&\leq S(r,f).
\end{align}
(3.1) subtract (3.8),
\begin{eqnarray}
(e^{p}-\psi)f+(\Delta_{\eta}^{n}a)^{(k)}-a\psi\equiv0.
\end{eqnarray}
We discuss following two cases.\\

{\bf Case 1} \quad $e^{p}\not\equiv\psi$. Then by (3.12) we obtain $T(r,f)=S(r,f)$, a contradiction.\\

{\bf Case 2} \quad $e^{p}\equiv\psi$. Easy to see from (3.12) that
 \begin{eqnarray}
(\Delta_{\eta}^{n}a)^{(k)}=e^{p}a.
\end{eqnarray}
We claim that $p$ is a constant. Otherwise,  if $p$ is a nonconstant polynomial, then $\lambda(e^{p})\geq1$. It follows from Lemma 2.4 and   that $a$ is a rational function that
 \begin{eqnarray}
\lambda(e^{p})=\lambda(\frac{(\Delta_{\eta}^{n}a)^{(k)}}{a})\leq max\{a,(\Delta_{\eta}^{n}a)^{(k)}\}<1,
\end{eqnarray}
which contradicts with $\lambda(e^{p})\geq1$. Thus, $p$ is a constant. We set $e^{p}=c$ to be a nonzero constant, then (3.13) implies $(\Delta_{\eta}^{n}a)^{(k)}=ca$. However, since $a$ is a rational function, and according to a simple calculation, we can not obtain $(\Delta_{\eta}^{n}a)^{(k)}=ca$, a contradiction.

From above discussion, we have $f\equiv(\Delta_{\eta}^{n}f)^{(k)}$. This completes the proof of Theorem 1.

\

{\bf Acknowledgements} The author would like to thank to anonymous referees for their helpful comments.



\begin{thebibliography}{99}
\bibitem{cx1} S. J. Chen, A. Z. Xu, \emph{Uniqueness theorem for a meromorphic function and its exact difference}, Bull. Korean Math. Soc. 57 (2020), no. 5, 1307-1317.
\bibitem{cy}  Z. X. Chen, H. X. Yi,  \emph{On Sharing Values of Meromorphic Functions and Their Differences},  Res. Math. 63 (2013),  557-565.

\bibitem{cf1}  Y. M. Chiang, S. J. Feng,
\emph{ On the Nevanlinna characteristic of $f(z+\eta) $ and
difference equations in the complex plane}, Ramanujan J. 16 (2008), no. 1,
105-129.

\bibitem{cf2} Y. M. Chiang, S. J. Feng, \emph{ On the growth of logarithmic differences, difference
quotients and logarithmic derivatives of meromorphic functions},
Trans. Amer. Math. Soc. 361 (2009), 3767-3791.

\bibitem{cc} N. Cui, Z. X. Chen, \emph{ The conjecture on unity of meromorphic functions concerning their differences}, J. Difference Equ. Appl. 22 (2016), no. 10, 1452šC1471.

\bibitem{h1} R. G. Halburd, R. J. Korhonen,
\emph{ Difference analogue of the lemma on the logaritheoremic derivative with
applications to difference equations}, J. Math. Anal. Appl. 314 (2006),
no. 2, 477-487.

\bibitem{h2} R. G. Halburd, R. J. Korhonen,
\emph{ Nevanlinna theory for the difference operator},
Ann. Acad. Sci. Fenn. Math. 31 (2006), no. 2, 463-478.

\bibitem{h3} R. G. Halburd, R. J. Korhonen and K. Tohge, \emph{ Holomorphic curves with shift-invarant hyperplane preimages},
Trans. Am. Math. Soc. 366 (2014), no. 8, 4267-4298.

\bibitem{h4} W. K. Hayman,\emph{ Meromorphic functions}, Oxford Mathematical Monographs Clarendon Press, Oxford 1964.

\bibitem{hkl} J. Heittokangas, R. Korhonen,  I. Laine, J. Rieppo, \emph{ Uniqueness of meromorphic functions sharing values with their shifts}, Complex Var. Elliptic Equ. 56 (2011), 81-92.

\bibitem{lyk} X. M. Li, H. X. Yi, C. Y. Kang \emph{ meromorphic functions sharing three values with their difference operators}, Bull. Korean Math. Soc. 52 (2015), no. 5, 1401šC1422.

\bibitem{ld1} K. Liu, X. J. Dong, \emph{ Some results related to complex differential-difference equations of certain types},
 Bull. Korean. Math. Soc. 51 (2014), 1453-1467.

\bibitem{ll} F. L\"u, W. R. L\"u, \emph{ meromorphic functions sharing three values with their difference operators}, Comput. Methods Funct. Theory 17 (2017), no. 3,  395-403.

\bibitem{qxg} X. G. Qi, \emph{ Value distribution and uniqueness of difference polynomials and entire solutions of difference equations},
 Ann. Polon. Math. 102 (2011), 129-142.

\bibitem{qly} X. G. Qi, N. Li, L. Z. Yang, \emph{ Uniqueness of meromorphic functions concerning their differences and solutions of difference PainlevšŠ equations}, Comput. Methods Funct. Theory 18 (2018), 567-582.

\bibitem{qy} X. G. Qi, L. Z. Yang, \emph{ Uniqueness of meromorphic functions concerning their shifts and derivatives}, Comput. Methods Funct. Theory 20 (2020), no. 1, 159-178.

\bibitem{y1} C. C. Yang, H. X. Yi,
\emph{ Uniqueness theory of meromorphic functions}, Kluwer Academic Publishers Group, Dordrecht, 2003.

\bibitem{y2} L. Yang, \emph{ Value Distribution Theory}, Springer-Verlag, Berlin, 1993.




\end{thebibliography}
\end{document}